\documentclass[12pt]{amsart}
\usepackage{amsfonts,graphicx,amssymb,amscd,amsmath,enumerate,verbatim,calc}

\def\NZQ{\mathbb}               

\def\QQ{{\NZQ Q}}
\def\ZZ{{\NZQ Z}}
\def\RR{{\NZQ R}}
\def\CC{{\NZQ C}}

\def\eb{{\bold e}}

\def\opn#1#2{\def#1{\operatorname{#2}}} 
\opn\pd{pd} 
\opn\rk{rk}
\opn\rank{rank}
\opn\depth{depth} 
\opn\grade{grade} 
\opn\height{height}
\opn\embdim{emb\,dim} 
\opn\codim{codim}
\opn\Tr{Tr} 
\opn\bigrank{big\,rank}
\opn\lcm{lcm}
\opn\reg{reg} 
\opn\ini{in} 
\opn\size{size}
\opn\mult{mult}
\opn\dist{dist}
\opn\cone{cone}
\opn\lex{lex}
\opn\rev{rev}
\opn\div{div}
\opn\Div{Div}
\opn\cl{cl}
\opn\Cl{Cl}
\opn\Syz{Syz} \opn\Im{Im} \opn\Ker{Ker} \opn\Coker{Coker}
\opn\Hom{Hom} \opn\Tor{Tor} \opn\Ext{Ext}
\opn\End{End} \opn\Aut{Aut} \opn\id{id} \opn\nat{nat}
\opn\mod{mod} \opn\ord{ord}
\opn\aff{aff} \opn\con{conv} \opn\relint{relint} \opn\st{st}
\opn\lk{lk} \opn\cn{cn} \opn\core{core} \opn\vol{vol}
\opn\link{link} \opn\star{star}

\def\Cc{{\mathcal C}}

\def\Fc{{\mathcal F}}
\def\Pc{{\mathcal P}}

\newtheorem{Theorem}{Theorem}[section]
\newtheorem{Lemma}[Theorem]{Lemma}

\newtheorem{Remark}[Theorem]{Remark}

\newtheorem{Example}[Theorem]{Example}

\newtheorem{Conjecture}[Theorem]{Conjecture}

\let\epsilon\varepsilon
\let\phi=\varphi
\let\kappa=\varkappa
\def\qed{\ifhmode\textqed\fi
      \ifmmode\ifinner\quad\qedsymbol\else\dispqed\fi\fi}
\def\textqed{\unskip\nobreak\penalty50
       \hskip2em\hbox{}\nobreak\hfil\qedsymbol
       \parfillskip=0pt \finalhyphendemerits=0}
\def\dispqed{\rlap{\qquad\qedsymbol}}

\begin{document}
\title{
Counterexamples of the conjecture \\
on roots of Ehrhart polynomials
}
\author{Akihiro Higashitani}
\thanks{
\;\;\;\; 
The author is supported by JSPS Research Fellowship for Young Scientists. 
}
\address{Akihiro Higashitani,
Department of Pure and Applied Mathematics,
Graduate School of Information Science and Technology,
Osaka University,
Toyonaka, Osaka 560-0043, Japan}
\email{a-higashitani@cr.math.sci.osaka-u.ac.jp}

\begin{abstract}
An outstanding conjecture on roots of Ehrhart polynomials says that 
all roots $\alpha$ of the Ehrhart polynomial of an integral convex polytope of dimension $d$ 
satisfy $-d \leq \Re(\alpha) \leq d-1$. 
In this paper, we suggest some counterexamples of this conjecture. 
\end{abstract}

\maketitle

\section*{Introduction}

Recently, on many papers, e.g., 
\cite{BDDPS}, \cite{BHW}, \cite{Bra}, \cite{BD} and \cite{MHNOH}, 
the root distributions of the Ehrhart polynomials 
have been studied intensively. 
In particular, one of the most significant problems is 
to solve the conjecture given in \cite[Conjecture 1.4]{BDDPS}. 
However, it will turn out that this conjecture is not true. 

First of all, we review what the Ehrhart polynomial is. 
Let $\Pc \subset \RR^N$ be an integral convex polytope 
of dimension $d$ and $\partial \Pc$ its boundary. 
Here an integral convex polytope is a convex polytope 
all of whose vertices have integer coordinates. 
Given positive integers $n$, we write 
$$i(\Pc, n) = |n \Pc \cap \ZZ^N|, \;\;\;\;\;\; 
i^*(\Pc, n) = |n (\Pc \setminus \partial\Pc) \cap \ZZ^N|,$$ 
where $n \Pc = \{ n \alpha : \alpha \in \Pc \}$ and 
$|X|$ denotes the cardinality of a finite set $X$. 
The systematic studies of $i(\Pc,n)$ originated in the work of Ehrhart \cite{Ehr}, 
who established the following fundamental properties: 
\begin{itemize}
\item $i(\Pc,n)$ is a polynomial in $n$ of degree $d$. 
(Thus, in particular, $i(\Pc,n)$ can be defined for every integer $n$, 
more generally, for every complex number $n$.) 
\item $i(\Pc,0) = 1$. 
\item (loi de r\'eciprocit\'e) 
$i^*(\Pc,n)=( - 1 )^d i(\Pc, - n)$ for every integer $n > 0$. 
\end{itemize}
We call this polynomial $i(\Pc,n)$ the {\em Ehrhart polynomial} of $\Pc$. 
We refer the reader to \cite[Part II]{HibiRedBook} and \cite[pp. 235--241]{StanleyEC} 
for the introduction to the theory of Ehrhart polynomials. 

We define the sequence $\delta_0, \delta_1, \delta_2, \ldots $ of integers 
by the formula 
\begin{eqnarray*}
(1 - \lambda)^{d + 1} \sum_{n=0}^{\infty} i(\Pc,n) \lambda^n 
= \sum_{j=0}^{\infty} \delta_j \lambda^j. 
\end{eqnarray*}
Since $i(\Pc, n)$ is a polynomial in $n$ of degree $d$ with $i(\Pc, 0) = 1$, 
a fundamental fact on generating functions 
(\cite[Corollary 4.3.1]{StanleyEC}) guarantees that $\delta_j = 0$ for every $j > d$. 
The sequence $\delta(\Pc) = (\delta_0, \delta_1, \ldots, \delta_d)$ 
is called the {\em $\delta$-vector} of $\Pc$. 
By the reciprocity law, one has 
\begin{eqnarray*}
\sum_{n=1}^{\infty} i^*(\Pc,n) \lambda^n = 
\frac{\sum_{i=0}^d \delta_{d-i} \lambda^{i+1}}{(1 - \lambda)^{d + 1}}. 
\end{eqnarray*}
The following properties on $\delta$-vectors are well known: 
\begin{itemize}
\item $\delta_0=1$ and $\delta_1 = |\Pc \cap \ZZ^N| - (d + 1)$. 
\item $\delta_d = |(\Pc \setminus \partial \Pc) \cap \ZZ^N|.$ 
Hence, we have $\delta_1 \geq \delta_d$. 
\item Each $\delta_i$ is nonnegative (\cite{StanleyDRCP}). 
\item When $d = N$, the leading coefficient $(\sum_{i=0}^d\delta_i)/d!$ of $i(\Pc,n)$ 
is equal to the usual volume of $\Pc$ (\cite[Proposition 4.6.30]{StanleyEC}). 
In general, the positive integer $\vol(\Pc) = \sum_{i=0}^d\delta_i$ 
is said to be the {\em normalized volume} of $\Pc$. 
\end{itemize}


For a complex number $a \in \CC$, let $\Re(a)$ denote the real part of $a$. 
Beck, De Loera, Develin, Pfeifle and Stanley propose the following 
\begin{Conjecture}\label{conj}{\em (\cite[Conjecture 1.4]{BDDPS})} 
All roots $\alpha$ of the Ehrhart polynomial of an integral convex polytope 
of dimension $d$ satisfy 
\begin{eqnarray}\label{ineq}
-d \leq \Re(\alpha) \leq d-1. \end{eqnarray}
\end{Conjecture}

It is proved in \cite{BDDPS} that this conjecture is true 
when $d=2$ and when roots are real numbers and it is also proved in \cite{BD} that 
this is also true when $d=3,4$ and 5. Moreover, in \cite{Bra}, 
the norm bound of roots of the Ehrhart polynomial is given with $O(d^2)$. 
In \cite{MHNOH}, for observing that this conjecture seems to be true, 
roots of the Ehrhart polynomials of several integral convex polytopes 
arising from finite graphs are discussed by using the languages of graph theory. 

In this paper, we show that Conjecture \ref{conj} is not true. (See Example \ref{ex}.) 
We can obtain many possible counterexamples by Theorem \ref{main} 
and we can find them for the first time when $d=15$.

\section{An important family of integral simplices}

This section is devoted to proving the following 
\begin{Theorem}\label{main}
Let $m,d,k \in \ZZ_{> 0}$ be arbitrary positive integers satisfying 
\begin{eqnarray}\label{condi}
m \geq 1, \; d \geq 2 \; \text{ and } \; 
1 \leq k \leq \lfloor (d+1)/2 \rfloor. 
\end{eqnarray}
Then there exists an integral convex polytope 
whose Ehrhart polynomial coincides with 
\begin{eqnarray}\label{ehrhart}
\binom{d+n}{d} + m \binom{d+n-k}{d}. 
\end{eqnarray}
\end{Theorem}

\subsection{How to compute the $\delta$-vector} 

Before proving the theorem, 
we recall from \cite[Part II]{HibiRedBook} 
the well-known combinatorial technique how to compute 
the $\delta$-vector of an integral simplex.

Given an integral $d$-simplex $\Fc \subset \RR ^N$ 
with the vertices $v_0, v_1, \ldots, v_d$, we set 
$$\widetilde \Fc=\left\{(\alpha,1)\in \RR^{N+1} \, : \, \alpha \in \Fc \right\},$$ 
which is an integral $d$-simplex in $\RR^{N+1}$ with the vertices 
$(v_0,1), (v_1,1), \ldots,(v_d,1)$. 
Clearly, we have $i(\Fc,n)=i(\widetilde \Fc,n)$ for all $n$. 
Let $$\Cc(\widetilde \Fc)=\Cc= 
\{r \beta : \beta \in \widetilde \Fc,0 \leq r \in \QQ \}.$$ 
Then one has 
$$i(\Fc,n)= \left| \left\{ (\alpha,n) \in \Cc : \alpha \in \ZZ^N \right\} \right|.$$ 
Each rational point $\alpha \in \Cc$ has a unique expression of the form 
$\alpha= \sum_{i=0}^{d}r_i(v_i,1)$ with each $0 \leq r_i \in \QQ$. 
Let $S$ be the set of all points $\alpha \in \Cc \cap \ZZ^{N+1}$ of the form 
$\alpha = \sum_{i=0}^{d}r_i(v_i,1),$ 
where each $r_i \in \QQ$ with $0 \leq r_i<1$. 
We define the degree of an integer point $(\alpha,n) \in \Cc$ with $\deg(\alpha,n)=n.$ 

\begin{Lemma}\label{compute}
Let $\delta_i$ be the number of integer points $\alpha \in S$ with $\deg \alpha=i$. Then 
$$\sum_{n=0}^{\infty}i(\Fc,n) \lambda^n 
=\frac{\delta_0+\delta_1\lambda+\cdots+\delta_d\lambda^d}{(1-\lambda)^{d+1}}.$$ 
\end{Lemma}

\subsection{A proof of Theorem \ref{main}}

We also recall the following well-known 
\begin{Lemma}\label{lemma}
Suppose that $(\delta_0,\delta_1,\ldots,\delta_d)$ is the $\delta$-vector 
of an integral convex polytope of dimension $d$. 
Then there exists an integral convex polytope of dimension $d+1$ 
whose $\delta$-vector is $(\delta_0,\delta_1,\ldots,\delta_d,0).$ 
\end{Lemma}

Now, we come to prove Theorem \ref{main}. Since we have 
$$
\sum_{n=0}^{\infty}\left( \binom{d+n}{d} + m \binom{d+n-k}{d} \right)\lambda^n 
= \frac{1 + m \lambda^k}{(1- \lambda)^{d+1}}, 
$$
it is sufficient to show that 
there exists an integral convex polytope $\Pc$ of dimension $d$ 
whose $\delta$-vector coincides with 
\begin{eqnarray*}
\delta_i=
\begin{cases}
1, \;\;\;\; &i=0, \\
m, &i=k, \\
0, &otherwise. 
\end{cases}
\end{eqnarray*}

When $k=1$, it is obvious that $(1,m,0,\ldots,0)$ is a possible $\delta$-vector. 
Thus, we assume that $k \geq 2$. In addition, by virtue of Lemma \ref{lemma}, 
our work is to find an integral convex polytope $\Pc$ of dimension $d$ 
with its $\delta$-vector 
\begin{eqnarray*}
\delta_i=
\begin{cases}
1, \;\;\;\; &i=0, \\
m, &i=(d+1)/2, \\
0, &otherwise, 
\end{cases}
\end{eqnarray*}
for arbitrary integers $m$ and $d$, where $m \geq 1$ and 
$d$ is an odd number with $d \geq 3$. 

Let $d \geq 3$ be an odd number and $c=(d-1)/2$. 
We define the integral $d$-simplex $\Pc \subset \RR^d$ 
by setting the convex hull of $v_0,v_1,\ldots,v_d$, which are of the form: 
\begin{eqnarray*}
v_i=
\begin{cases}
\eb_i, \quad\quad\quad\quad\quad\quad\quad\quad\quad\quad\quad 
&i=1,\ldots,d-1, \\
\sum_{j=1}^c\eb_j + \sum_{j=c+1}^{2c}m\eb_j + (m+1)\eb_d, 
&i=d, \\
(0,0,\ldots,0), &i=0, 
\end{cases}
\end{eqnarray*}
where $\eb_1,\eb_2,\ldots,\eb_d$ denote the unit coordinate vectors of $\RR^d$. 
In other words, for $i=1,2,\ldots,d$, $v_i$ is equal to the $i$th row vector of 
the $d \times d$ lower triangular integer matrix 
\begin{eqnarray}\label{mat}
\begin{pmatrix}
1      &0      &\cdots &\cdots &\cdots &\cdots &0 \\
0      &1      &\ddots &       &       &       &\vdots \\
\vdots &\ddots &\ddots &\ddots &       &       &\vdots \\
\vdots &       &\ddots &\ddots &\ddots &       &\vdots \\
\vdots &       &       &\ddots &\ddots &\ddots &\vdots \\
0      &\cdots &\cdots &\cdots &0      &1      &0 \\
1      &\cdots &1      &m      &\cdots &m      &m+1
\end{pmatrix}, 
\end{eqnarray}
where there are $c$ 1's and $c$ $m$'s in the $d$th row. 
Then we notice that $\vol(\Pc)=m+1$, 
which coincides with the determinant of \eqref{mat}. 

For $j=1,2,\ldots,m$, since 
$$\sum_{i=0}^c \frac{m+1-j}{m+1} (v_i,1) + \sum_{i=c+1}^d \frac{j}{m+1} (v_i,1) 
=(\underbrace{1,1,\ldots,1}_c,\underbrace{j,j,\ldots,j}_{c+1},c+1) \in \ZZ^{d+1}$$ 
and $$0 \leq \frac{m+1-j}{m+1}< 1, \;\;\; 0 \leq \frac{j}{m+1}<1, $$ 
Lemma \ref{compute} guarantees that $\delta_{c+1} \geq m$. 
Moreover, thanks to $\vol(\Pc)=m+1$ together with the nonnegativity of $\delta$-vectors, 
we obtain $\delta_{(d+1)/2}=m$. 
Therefore, we can conclude that $\Pc$ has the required $\delta$-vector.

\section{Counterexamples of Conjecture \ref{conj}}

In this section, we consider the roots of the polynomial \eqref{ehrhart} 
with positive integers $m,d$ and $k$ satisfying \eqref{condi}. 

Let $f(n)$ be the polynomial \eqref{ehrhart} in $n$ of degree $d$. 
Since 
$$f(n)= \frac{\prod_{j=1}^{d-k}(n+j)}{d!}
\left( \prod_{j=d-k+1}^d (n+j) + m \prod_{j=0}^{k-1}(n-j) \right), 
$$
$-1,-2,\ldots,-d+k$ are always the roots of $f(n)$. 
Hence, we consider the roots of $g_{m,d,k}(n)$ 
with positive integers $m,d$ and $k$ satisfying \eqref{condi}, 
where $g_{m,d,k}(n)$ is the polynomial 
$$g_{m,d,k}(n)=\prod_{j=d-k+1}^d (n+j) + m \prod_{j=0}^{k-1}(n-j)$$ 
in $n$ of degree $k$. 

\begin{Example}\label{ex}{\em 
Let us consider the polynomial $g_{m,15,8}(n)$. 
When $1 \leq m \leq 8$, all their roots satisfy \eqref{ineq}. 
On the contrary, when $m=9$, its eight roots are approximately 
\begin{eqnarray*}
&&14.37537447 \pm 25.02096544 \sqrt{-1}, \;\;\; 
-0.77681486 \pm 10.23552765 \sqrt{-1}, \\
&&-2.56596317 \pm 4.52757516 \sqrt{-1} \;\;\text{and}\;\; 
-3.03259644 \pm 1.31223697 \sqrt{-1}. 
\end{eqnarray*}
By virtue of Theorem \ref{main}, 
this implies that there is a counterexample of Conjecture \ref{conj}. 
Moreover, in the similar way, it can be verified that 
for every $15 \leq d \leq 100$, 
there is a root of $g_{9,d,\lfloor (d+1)/2 \rfloor}(n)$ 
which does not satisfy \eqref{ineq}, i.e., 
there is a counterexample of Conjecture \ref{conj} 
for each dimension $d$ with $15 \leq d \leq 100$. 
(Those are computed by {\tt Maple}.) 
It also seems to be true when $d \geq 101$. 
In addition, when $d \geq 17$, we can also verify that 
there is a root of $g_{9,d,\lfloor (d+1)/2 \rfloor}(n)$ 
whose real part is greater than $d$. 
}\end{Example}

Now, these computational results are also supported theoretically. 
In fact, for example on the roots of $g_{9,15,8}(n)$, 
by applying the {\em Routh-Hurwitz stability criterion}, 
we can check that there exists a root of $g_{9,15,8}(n+14.3)$ 
whose real part is nonnegative and 
the real parts of the roots of $g_{9,15,8}(n+14.4)$ are all negative. 
Of course, this means that there exists a root $\alpha$ of $g_{9,15,8}(n)$ 
with $14.3 \leq \Re(\alpha) < 14.4$.

\begin{Remark}{\em 
On the order of the largest real part of 
the roots of $g_{9,d,\lfloor (d+1)/2 \rfloor}(n)$, 
the order does not seem to be linear in $d$. 
For example, the largest real part is around 59 when $d=30$, 
it is around 174 when $d=50$ and it is around 722 when $d=100$. 
Thus, it is natural to claim that 
the upper bound of the real parts of the roots of the Ehrhart polynomials of 
integral convex polytopes is not $d-1$ but something with $O(d^2)$, 
while we do not know the lower bound. 
}\end{Remark}

\begin{Remark}{\em 
(a) When $m=1$, the real parts of all the roots of $g_{1,d,k}(n)$ 
coincide with $(-d+k-1)/2$, which satisfies $-d < (-d+k-1)/2 < -1/2$. 
In fact, since all the roots of $1+\lambda^k$ 
are on the unit circle, we can apply the theorem of \cite{RV} to 
the polynomial $\binom{n+d}{d}+\binom{n+d-k}{d}$. 
When $m=2$, on the other hand, 
we can obtain an other counterexample of Conjecture \ref{conj} 
when $d=37$ and $k=19$. \\
(b) When $k=1$, one has $g_{m,d,1}(n)=(m+1)n+d$. 
Thus, its root is $-d/(m+1)$, which satisfies $-d < -d/(m+1) <0$. 
When $k=2$, then one has $g_{m,d,2}(n)=(m+1)n^2+(2d-m-1)n+d(d-1)$. 
Let $D(g_{m,d,2}(n))$ denote the discriminant of $g_{m,d,2}(n)$. 
If $D(g_{m,d,2}(n))<0$, then the real part of the roots of $g_{m,d,2}(n)$ is 
$-d/(m+1)+1/2$, which satisfies $-d+1/2 < -d/(m+1)+1/2 < 1/2.$ 
Note that when we let $m$ grow sufficiently compared with $d$, 
the roots of $g_{m,d,2}(n)$ become real and they approach 0 and 1 respectively. 
In fact, the roots of $g_{m,d,2}(n)$ coincide with that of 
$g_{m,d,2}(n)/m=(n+d)(n+d-1)/m+n(n-1)$. 
}\end{Remark}

\smallskip

\section*{Acknowledgemenets} 
The author would like to thank Hidefumi Ohsugi and Tetsushi Matsui 
for giving him some comments on Example \ref{ex}, 
pointing out a gap between approximately roots and actual roots and 
telling him the criterion.

\end{document}